\documentstyle[12pt]{amsart}
\newtheorem{lem}{{\sc Lemme}}

\newtheorem{prop}{{\sc Proposition}}

\newtheorem{cor}{{\sc Corollaire}}

\newcommand{\rep}{repr\'esentation}

\newcommand{\2}{{\frac{1}{2}}}

\newcommand{\B}{{\bf B}}
\newcommand{\z}{\zeta}
\newcommand{\w}{\omega}
\newcommand{\s}{\langle}
\let\r\rangle

\begin{document}

\title
[repr\'esentations de Weil sur un corps fini]
{Construction g\'eom\'etrique de repr\'esentations de Weil sur un corps fini}

\author{Jos\'e Pantoja et Jorge Soto-Andrade}
\address{Pantoja: Instituto de Matem\'atica\\
Universidad Cat\'olica de Valpara\'{\i}so\\
Casilla 4059\\
Valpara\'{\i}so, Chile}
\email{jpantoja@@aix1.ucv.cl}
\address{Soto-Andrade: Departamento de Matem\'aticas\\
Universidad de Chile\\
Casilla 653\\
Santiago, Chile}
\email{sotoandr@@orphee.polytechnique.fr ou
sotoandr@@abello.dic.uchile.cl}
\thanks{Pantoja a \'et\'e subventionn\'e par DGI-UCV 124.736/94, par la
Coop\'eration fran\c caise (Programme ECOS), et par FONDECYT (Projet
1950895).
Soto-Andrade a \'et\'e subventionn\'e par la 
Coop\'eration fran\c caise (Programme ECOS), par FONDECYT (Projet
1940590), ainsi que par le Grant NSF DMS/9022140 pendant son s\'ejour au
 MSRI.}

\maketitle

\begin{abstract}
We construct, by contraction of a suitable complex vector bundle, the Weil
representation of the finite symplectic group $Sp(A)$. We give an explicit
description of the space of all lagrangian subspaces, which we use  to
compute the cocycle of our representation in terms of a geometric Gauss
sum. We recover in this way previously constructed generalized Weil
representations (see \cite{ast,cor}) by
restriction of our representation to an appropiate embedding of $SL(n) $
into $Sp(A)$.
\end{abstract}

\section{Connexions \'equivariantes sur un $Sp(A) $-fibr\'e vectoriel}
 
\indent
 
Soit $ (W,A) $ un espace   symplectique non-d\'eg\'en\'er\'e, de dimension paire  $n = 2m$ sur 
un corps fini  $k$  \`a  $q$  \'el\'ements. 
	Nous rappelons la construction g\'eom\'etrique de la \rep  de Weil
 du groupe  $Sp(A)$, par 
contraction d'un  $Sp(A)$ - fibr\'e vectoriel complexe convenable \`a l'aide d'une connexion
\'equivariante, qui a \'et\'e  indiqu\'ee en \cite{arc} dans le cas $m = 1$.

	Posons  $ G = Sp(A)$.
 D\'esignons $ \cal X$ l'ensemble de tous les caract\`eres non nuls du
 groupe additif $k^+$ de  $k$
et fixons un caract\`ere non-trivial $\psi \in \cal X$.  Nous
 d\'efinissons un  $G$ - fibr\'e vectoriel complexe  $(E,p,B, \tau)$
 par les donn\'ees suivantes:

- la base  $B$ est l'ensemble de tous les sous-espaces lagrangiens (c'est-\`a-dire totalement 
isotropes maximaux) de l'espace symplectique $W$;

- l'espace total  $ E$ est la r\'eunion disjointe des espaces $E_L  \; (L \in B) $ form\'es de 
toutes les fonctions complexes $f$ sur $W$ telles que 
$$
f(x+y ) = \psi(A(x,y))f(x )
$$
pour tous les $x \in W$  et  $y \in L, \psi \in  \cal X $ \'etant fix\'e;

- $p$ d\'esigne la projection canonique de $E$ sur $B$, qui  \`a chaque
$f \in E_L$ fait correspondre $L$;

- L'action  $\tau $  de $G$ dans $E$ et $B$ est donn\'ee par
$$  [\tau_g(f)](x) = f(g^{-1}(x)) $$,
pour  $g \in G, f \in E,  x \in  W,$    et
 $$  \tau_g(L) = g(L) $$
pour  $g \in G, L\in B $.
Nous d\'efinissons une connexion  $\Gamma = \{ \gamma_{L',L}\}_{L',L \in B}$ sur le fibr\'e 
 $(E,p,B) $ par
 
$$(\gamma_{L',L}f)(\omega,\psi) = (|L||L\cap L'|)^{-\2}\sum_{ \zeta \in L'}\psi(A(\omega, \zeta))f(\omega + \zeta)$$

%\newpage
\begin{prop}
{\bf Propri\'et\'es de $\Gamma$}
\begin{description}
\item[i)]
$$\gamma_{L,L} = Id_{E_L}$$
\item[ii)]
$$ \gamma_{L,L' }\circ\gamma_{L',L } = \gamma_{L,L}= Id$$ 
\item[iii)]
$$  \gamma_{L'',L' } \circ\gamma_{L',L } = |L'' \cap L'|^\2
(|L\cap L''|  |L' \cap L | |L | )^{- \2}{\cal S} _L(L',L'')\gamma_{L'',L } $$
o\`u l'on note ${\cal S} _L(L',L'')$ , ou plus pr\'ecisement ${\cal S}^W _L(L',L'')$, la somme de
Gauss g\'eom\'etrique associ\'ee \`a l'espace symplectique $W$ et ses lagrangiens $L, L', L'' $ ,
d\'efinie par 

%$$ \stackrel{^W}{ }\hspace{-7pt}{\cal S} _L(L',L'') = {\cal S} _L(L',L'') =  
%\sum _{\zeta \in L \cap
%(L' + L'')} \psi(A(\zeta', \zeta'')) \; ,$$ 

$${\cal S} _L(L',L'') = {\cal S}^W _L(L',L'') = \sum _{\zeta \in L \cap
(L' + L'')} \psi(A(\zeta', \zeta'')) \; ,$$  

%\sum _{\zeta \in L \cap
%(L' + L'')} \psi(A(\zeta', \zeta'')) \; ,$$ 
o\`u $\zeta'$ (resp. $\zeta''$ ) $\;$ d\'esigne la
composante de $\zeta  \in  L \cap (L' + L'')$. selon $L'$ (resp. $L''$).
\end{description} 
\end{prop}
\begin{flushleft}

{\em D\'emostration.-} 
\end{flushleft}

ii) On a

$$(\gamma_{L,L' }\circ\gamma_{L',L }f)(\w) = (|L||L\cap L'|)^{-\2} \sum_{ \zeta \in L}\psi(-A(\omega,
\zeta))(\gamma_{L',L }f)(\w+ \z)  = $$
$$ = (|L||L\cap L'|)^{-1}\sum_{ \zeta \in L} \sum_{ \zeta \in L'}\psi(-A(\omega,
\zeta) - A(\w+\z,\z'))( f)(\w+ \z + \z')$$
$$= (|L||L\cap L'|)^{-1} \sum_{ \zeta \in L'}\sum_{ \zeta \in L}\psi(-A(\omega,
\zeta) - A(\w+\z,\z')+ A(\w+\z',\z))( f)(\w + \z')$$
$$
= (|L||L\cap L'|)^{-1} \sum_{ \zeta \in L'}[\sum_{ \zeta \in L}\psi(2A(\z',\z)]\psi( 
  - A(\w,\z') ) f(\w + \z')
  $$
$$
  = (|L\cap L'|)^{-1} \sum_{ \zeta \in L'\cap L} \psi( 
  - A(\w,\z') ) f(\w + \z') = f(\w) \; , $$
  puisque le charact\`ere $\z \mapsto  \psi(A(\z',\z)$  du groupe additif de $L$ est non-trivial si et
seulement si  $\z' \in L' \cap L^\perp = L \cap L'$ et que $f \in E_L$

iii)
Nous avons, quels que soient $f \in E_L, \w \in W , \psi \in \cal X$, 

\begin{eqnarray*}
\lefteqn{ (\gamma_{L,L'' }\circ\gamma_{L'',L' } \circ\gamma_{L',L }f)(\w )   =}\\
 &= & (|L|^3|L''\cap L||L'\cap L''|)|L\cap L'|)^{-\2}\sum_{   \zeta'' \in L''}  \sum_ {\z' \in L'
} \sum_{\z \in L}  \psi( 
  - A(\w,\z  ) -A(\w+\z, \z'') - \\
 &  & - A(\w + \z + \z'', \z'))f(\w +\z + \z'' +\z')\\
  &  = & (|L|^3|L''\cap L||L'\cap L''||L\cap L'|)^{-\2}\sum_{   \zeta'' \in L''}  \sum_ {\z' \in L'
} [\sum_{\z \in L}  \psi(2A(\z''+\z', \z))]\psi( - A(\w,\z' + \z''  ) -A(\z'', \z')) \\
& & f(\w + \z'' +\z') \\
& = & (|L||L''\cap L||L'\cap L''||L\cap L'|)^{-\2} [\sum_{  \begin{array}{c}  \z' \in L',
\z'' \in L'' \\  
\z' + \z'' \in  L
\end{array} }  \psi(A(\z', \z''))     ] \psi( - A(\w,\z' + \z''  )  f(\w + \z'' +\z') \\
& = & (|L||L''\cap L||L\cap L'|)^{-\2} |L'\cap L''|^{\2}
[\sum_{    \z \in L \cap (L' + L'')}  
\psi(A(\z', \z''))]       
\end{eqnarray*} 
o\`u, pour chaque $\z \in L \cap (L' + L'') $, on note  $  \z'$(resp. $ \z''$) sa composante dans
$L'$ (resp.$ L''$) selon  une d\'ecomposition quelconque  $\z = \z' + \z'' .$

On d\'efinit ainsi une forme quadratique $ Q_L^{L',L''}$ sur l'espace $L \cap (L' + L'') $ par
 $$Q_L^{L',L''}(\z) = A(\z , \z') , $$
 quel que soit  $\z \in   L \cap (L' + L'') $ ' , et  o\`u  $\z = \z' + \z'' $ est une d\'ecomposition
quelconque de  $\z \in  L \cap (L' + L'') $, avec  $\z' \in L'$ et  $\z' \in L''$. Notons que la
valeur de  $A(\z' , \z'')  $ ne d\'epend alors pas de la d\'ecomposition choisie de $\z \in L \cap
(L' + L'')$

\begin{cor}.- $\Gamma$ est une connexion $G$ -\'equivariante sur le fibr\'e $ (E,p,B, \tau)$
, dont le multiplicateur $\mu$ est donn\'e par
$$ \mu(L'',L",L) = |L'' \cap L'|^\2
(|L\cap L''|  |L' \cap L | |L | )^{- \2}{\cal S} _L(L',L'')
$$
\end{cor}
\section { Calcul de la somme de Gauss g\'eom\'etrique ${\cal S} _L(L',L'')$}

\indent

	Nous d\'ecrivons tout d'abord les sous-espaces lagrangiens de l'espace symplectique $(W,B)$.
% \section{Modules symplectiques.-}

 Notons  $A$ l'anneau 
involutif $A$ des matrices $n \times  n$ sur  $k$, 
dont l'involution est la 
transpos\'ee. A chaque     espace symplectique $(W,B)$, de dimension paire $2m$, sur le corps $k$
on associe  un 
 module symplectique $(M, \B)$ sur $A$, comme suit:\\
On pose   $M = W^m$ et l'on muni  $M$ de la forme sesquilin\'eaire altern\'ee $ \bf B$ \`a valeurs
dans $A$ donn\'ee par 
$$ {\bf B}(u,v) = (B(u_i,v_j))_{1\leq i,j \leq m} $$
On a alors
$$ \B(hu,v) = h\B(u,v) \; $$
$$ \B(u,hv) = \B(u,v)h^*   $$
$$  \B(v,u) = \B(u,v)^* \;$$
D'autre part d\'esignons par 
$\langle u \rangle$  le sous-espace  vectoriel de $W$  engendr\'e par les composantes $u_1, \dots ,u_n$
de $u \in M = W^m$. On note alors $rg u$ le rang de $u \in M$, d\'efini comme la dimension de l'espace
$\s u \r$. On a $ rg u = m$ si et seulement si $u$ est un vecteur libre dans le  $A$-module $M$. Notons
 que $\s u \r  = \s v \r$
 si et seulement si  $ v = hu $ pour une matrice  $h \in A$ convenable.
	Remarquons que pour que $\s u \r$ soit  Lagrangien il faut et il suffit que $u$ soit un vecteur
isotrope libre dans le module symplectique $(M, \B)$.
Soit $\{p_1, \dots, p_m, q_1, \dots , q_m \}$ une base symplectique de $(W,B)$. Alors on a 
$$ B(p_i,q_j) = \delta_{i,j}
$$
quels que soient $ 1 \leq i,j \leq m $ et les deux vecteurs isotropes    $P = (p_1, \dots, p_m)$ et $Q
=  (q_1, \dots , q_m)$ forment
 un syst\`eme libre dans  $M$ tel que $ \B(P,Q) = Id_m$  
 \begin{prop}
Pour que l'espace $ L_{a,b} = \s aP + bQ \r,  \; \; (a,b \in A)$ soit Lagrangien il faut et il
 suffit que  $aA + bA = A$ (on dit alors que $a$ et $b$ sont relativement premiers) et que 
 $ab^* = ba^*$. En plus, tout sous-espace Lagrangien  $L$ de $(W,B)$ est un   $L_{a,b} $
 pour $a, b \in A $  convenables et  $L_{a,b}  =  L_{a',b'}$  \'equivaut \`a  $ a' = ca  $ et $ b' = cb$ pour un $c \in A$.

\end{prop}
    
%Developper le point de vue anneau involutif et module symplectique sur
%un anneau involutif.
Remarquons que le $A$-module  $M$ est muni d'une application 
$ (\lambda, v) \mapsto  \lambda \cdot u = \lambda_1u_1 +  \dots + 
\lambda_mu_m $ de  $k^m \times M$ sur $W$, qui permet d'identifier en fait $M$ avec 
%le produit tensoriel  $k^m \otimes W$ de $k^m$ et $W$ sur $k$
$Hom_k(k^m, W)$ en tant que $A$- module gauche.

\begin{lem}

On a 
$$  g(L_{[a,b]}) = L_{[a,b]g^*}  $$
quelques soient  $ g \in G, L_{[a,b]}  \in \cal L$
\end{lem}
\begin{lem}
On a, quels que soient $\lambda, \mu \in k^m, a,b \in A, u,v \in M,$ 
$$ B(\lambda \cdot u,\; \;  \mu \cdot v )= \lambda \ \B(u,v) \mu $$
\end{lem}
 
\begin{prop}

On a
$$    S_{L_{[a,b]}}(L',L'') = S(ab^ {\ast}) , $$
o\`u l'on note $S(Q)$ la somme de Gauss classique 
$$  \sum_{x \in E} \psi(Q(x))  \; ,$$  associ\'ee \`a  un espace quadratique $(E,Q) \; ,$
le caract\`ere non-trivial $\;\psi$ du groupe additif $\;
k^+$ \'etant fix\'e.
\end{prop}
\begin{lem}
Soient $  L', L'' \in \cal L .\;\;$ 
Alors pour chaque $L \in \cal L$ le sous-espace 
$$ \tilde{L} = 
 L \cap (L' + L '') + 
 L' \cap L'' $$  de $L' + L ''$  
est encore un sous-espace lagrangien de $W$ et le sous-espace
$\tilde{L}/ L' \cap L''$ est un lagrangien de l'espace symplectique non-d\'eg\'en\'er\'e $(L' + L
'')/  L' \cap L''$. En plus on a

$$
S_{L}^W(L',L'') = |L \cap L' \cap L''| S^{(L' + L'')/  L' \cap L''}_{\tilde{L}/ L' \cap
L''}({L'/ L' \cap L'', L''/ L' \cap L''}) 
 $$ 
\end{lem}

\begin{prop}

On a
$$
|S_{L}(L',L'')| = |L'' \cap L'|^{-\2}(|L\cap L''|  |L' \cap L | |L | )^{ \2} = q^{m - \2r(ab^*)} =
q^{\2(3m - r(a) - r (b))}, 
$$
quels que soient $L, L', L'' \in \cal L $, o\`u $r(x)$ d\'esigne le rang d'une matrice $x \in A$.

\end{prop}
\begin{prop}
Le cocycle $c(g,h) $ de  la repr\'esentation de Weil  $(V,\rho)$ de $G$ associ\'ee \`a  $L_o \in B$
est donn\'e par
$$
c(g,h) =  |{\cal S} _{L_o}(gL_o,ghL_o)|^{-1}{\cal S} _{L_o}(gL_o,ghL_o)  \;
$$
quels que soient  $ g,h \in  G$.
\end{prop}

%Valeur  explicite du cocycle   $c(g,h) $  en termes de la somme de Gauss...

\section{Les repr\'esentations de Weil g\'en\'eralis\'ees des groupes 
$ SL(n,k)$} 

\indent

	Soit $V$ un espace vectoriel de dimension finie $n$ sur le corps fini $k$. Nous d\'efinissons 
l'espace symplectique $ (W, A) $ associ\'e \`a  $V$  par 
$ \;
W =  \bigwedge^1 V  \oplus   \bigwedge^3 V \oplus \dots  \oplus 
  \bigwedge^{n-1}V 
\;\;\;$
si  $n$ est pair (resp. $\;W =  \bigwedge^1 V \oplus   \bigwedge^2 V \oplus \dots  \oplus   \bigwedge^{n-1}V 
\;\;\;$ si  $n$  est impair) et  
$$   A(\omega, \zeta) = (\omega \wedge \zeta)_n \in   \bigwedge^n V \simeq k  $$
(resp. 

$$ A(\omega, \zeta) = (\omega \wedge \zeta^\vee)_n \in 
  \bigwedge^n V \simeq  k  \;\; , $$
o\`u  $ \zeta^\vee = \sum_{0 \leq i \leq n }(-1)^i\zeta_i \; $ pour   $ \;\;
\zeta = \sum_{0 \leq i \leq n }\zeta_i ).$

Le groupe  $H = SL(V) \simeq SL(n,k)$   agit naturellement  dans $(W,A)$ et s'identifie 
ainsi a un sous-groupe de $G$ , auquel nous pouvons restreindre la repr\'esentation de 
de Weil g\'en\'eralis\'ee d\'ecrite ci-dessus. Notons $(U, \rho)$ la repr\'esentation de $H$ ainsi
obtenue. 

On recup\`ere ainsi les repr\'esentations que nous avons
construites  dans \cite{ast}, \cite{cor}. 
De plus, comme la connexion $\Gamma$ que nous avons construite ci-dessus
 permet de passer d'une $H$
orbite quelconque \`a une autre, elle fournit en fait des isomorphismes
 entre des repr\'esentations 
de Weil g\'en\'eralis\'ees de $H$ associ\'ees \`a des orbites
diff\'erentes.  Ceci s'applique en 
particulier au cas des $H$-orbites r\'eduites \`a un lagrangien et aux
 orbites introduites dans loc.
cit., et fournit une construction uniforme de certains isomorphismes
 notables entre des repr\'esentations
naturelles et des repr\'esentations de Weil g\'en\'eralis\'ees,
 d\'ej\`a remarqu\'es par P. Cartier dans
les ann\'ees 60, dans le cas $n = 2 $.    

Nous remarquons enfin que l'on peut construire des op\'erateurs
 d'entrelacement de $(U, \rho)$  comme suit.

Soit  $\phi$ un automorphisme symplectique $W$ qui commute \`a
 l'action de $H$ dans $W$. L'automorphisme
$\phi$ fournit   alors un automorphisme du   $G$-fibr\'e vectoriel
 associ\'e \`a l'espace
symplectique $(W,A)$, d\'efini au paragraphe 1, qui commute avec sa
 connection  $\Gamma$. On en 
tire aussit\^ot un automorphisme  $\Phi$ de la repr\'esentation  $(U, \rho)$.

On peut v\'erifier que, pour $ n \leq 4 $, 
 le groupe d'op\'erateurs d'entrelacement ainsi obtenu (qui repr\'esente
donc le centralisateur $Z_G(H)$  de $H$  dans $G$ ) engendre 
l'essentiel de l'alg\`ebre commutante de 
la repr\'esentation $(U, \rho)$ de $H$ (au sens que le rapport entre
 la dimension du
sous-espace vectoriel engendr\'e par ces op\'erateurs et la dimension 
de l'alg\`ebre commutante de 
la repr\'esentation $(U, \rho)$ tend vers $1$ pour $q$ tendant vers
l'infini).

\end{document}